\documentclass{amsart}
\usepackage{amscd}
\usepackage{amssymb}
\theoremstyle{plain}
\newtheorem{theorem}{Theorem}[section]
\newtheorem{lemma}[theorem]{Lemma}
\newtheorem{proposition}[theorem]{Proposition}
\newtheorem{corollary}[theorem]{Corollary}

\theoremstyle{definition}
\newtheorem{definition}[theorem]{Definition}
\newtheorem{def-prop}[theorem]{Definition-Proposition}
\newtheorem{example}[theorem]{Example}

\theoremstyle{remark}
\newtheorem{remark}{Remark}[section]
\begin{document}
\title[Signed Quivers]{ Signed quivers, symmetric quivers,  \\
and root systems.  }
\email{mitia@mccme.ru}
\author{D.A.Shmelkin}
\begin{abstract}
We define a special sort of weighted oriented graphs, signed quivers.
Each of these yields a symmetric quiver, i.e., a quiver endowed with
an involutive anti-automorphism and the inherited signs.
We develop a representation theory of symmetric quivers, in particular
we describe the indecomposable symmetric representations.
Their dimensions constitute  root systems corresponding
to  certain symmetrizable generalized Cartan matrices.
\end{abstract}
\maketitle

\section{Introduction.}

Let $Q$ be a finite quiver, i.e., an oriented graph. 
Fix the notation as follows. We denote by $Q_0$ and $Q_1$ the sets
of the vertices and the arrows of $Q$, respectively. For any arrow
$\varphi\in Q_1$ denote by $t\varphi$ and $h\varphi$ its tail and
its head, respectively. A representation $V$ of $Q$ over ${\bf k}$ 
consists in defining  a vector space $V(i)$ over ${\bf k}$, 
for any $i\in Q_0$, and  a ${\bf k}$-linear map 
$V(\varphi):V(t\varphi)\to V(h\varphi)$, for any 
$\varphi\in Q_1$. The dimension vector $\dim V$ is the collection
of $\dim V(i), i\in Q_0$. For a fixed dimension $\alpha$, we may
fix $V(i)={\bf k}^{\alpha_i}$. Then the set $R(Q,\alpha)$ of 
the representations of dimension $\alpha$ is converted into the
vector space
\begin{equation}
R(Q,\alpha) = \bigoplus_{\varphi\in Q_1} {\rm Hom}
({\bf k}^{\alpha_{t\varphi}},{\bf k}^{\alpha_{h\varphi}}).
\end{equation}
A homomorphism $H$ of a representation $U$ of $Q$ to another representation,
$V$ is a collection of linear maps $H(i), U(i)\to V(i)\in Q_0$ such that
for any $\varphi\in Q_1$ holds $V(\varphi) H(t\varphi) = H(h\varphi)U(\varphi)$.
The endomorphisms, automorphisms, and isomorphisms are defined naturally.
An easy but very fruitful observation is that the isomorphism classes
of representations of $Q$ are the orbits of a reductive group
\begin{equation}
GL(\alpha) = \prod_{i\in Q_0} GL(\alpha_i)
\end{equation}
acting naturally on $R(Q,\alpha)$, as follows

\begin{equation}
(g(V))(\varphi) = g(h\varphi) V(\varphi) (g(t\varphi))^{-1}.
\end{equation} 
Futhermore, by the Krull-Schmidt theorem, each representation has a unique,
modulo isomorphisms and permutations of summands, decomposition into
{\it indecomposable} ones. The classification of representations
modulo isomporphism is therefore reduced to that for the indecomposable
ones. The latter problem is solved for the {\it finite} and {\it tame} quivers.
Moreover,  by Kac's Theorem \cite{kac2} the dimensions
of the indecomposable representations are exactly the positive roots
of the symmetric Kac-Moody algebra corresponding to the underlying graph of $Q$.

From the point of view of Invariant Theory, the above notation
introduces a nice set of reductive linear groups,
$(GL(\alpha),R(Q,\alpha))$. This set is nice because the developed language
of quiver representations allows (at least in some cases) to describe
orbits, invariants, semi-invariants etc. Another important feature of this
set is that  the underlying quiver $Q$
determines many natural properties of the groups, for example,
if $Q$ is finite, then  $(GL(\alpha),R(Q,\alpha))$ contains finitely many orbits for
any $\alpha$. 

The above set of groups was extended in \cite{dw1} to that of {\it generalized quivers}.
These can be described as actions of certain reductive groups in the spaces of 
either orthogonal or symplectic representations
of a {\it symmetric} quiver $S$, i.e. a quiver with an involutive anti-automorphism.
In this setting one can also classify the orbits in terms of
the indecomposable representations; moreover, an important result of \cite[Proposition~2.7]{dw1}
states that each indecomposable 
orthogonal or symplectic representation
is either indecomposable as a representation of $S$, or is a sum of two
indecomposable representations in involution. 

We suggest another extension of the set of groups containing all the generalized quivers. 
Similarly to the usual quiver context, these groups are parametrized by
a {\it signed quiver} and a dimension vector. We represent these groups
as actions on the spaces of {\it symmetric} representations of a symmetric  (in the same sense
as in \cite{dw1}) quiver endowed with signs. A key observation
\ref{isofiso} is that symmetric representations are conjugate by our group
if and only if they are isomorphic as quiver representations. 
We obtain an a priori description  \ref{1_or_2} of
the indecomposable symmetric representations in the same style as in \cite{dw1}.
As a corrolary we get in \ref{class_indec} the uniqueness of the decomposition for the symmetric representations.  

Analogously to the usual quiver context, the dimensions
of indecomposable symmetric representations turn out to be related to the root system
for a {\it symmetrizable} generalized Cartan matrix. 
Let $\Gamma$ be the underlying graph of our symmetric quiver and  
assume that the corresponding involutive automorphism $\pi$ is {\it admissible} (see Definition \ref{admis}).
Then $\pi$ acts on the root system $\Delta(\Gamma)$ and
one may consider the set
consisting of the $\pi$-invariant roots and
of the sums $\alpha+\pi(\alpha)$ for non-invariant roots $\alpha$.
Recently it was observed in \cite{hu} that this set is the 
root system of symmetrizable type corresponding to a graph, which is
a sort of factorization of $\Gamma$ by $\pi$. In fact this observation
goes back to Kac: we can not find it as a statement but it is used in
\cite[7.9]{kac3} in order to construct the root systems 
for $B_l,C_l,F_4,G_2$ in terms of those for $D_{l+1},A_{2l-1},E_6,D_4$, respectively. 
The statement from \cite{hu} concerns not all signed quivers because it
does not depend on the signs. We then generalize this: we define
a graph $\Gamma_{\pi,\sigma}$ depending on both
$\pi$ and the signs $\sigma$ and obtain in \ref{root_system} a description
of the root system $\Delta(\Gamma_{\pi,\sigma})$ in terms of $\Delta(\Gamma)$, 
similar to the above one. In \ref{dim_sym_rep} we give a sufficient condition
for the set of dimensions of indecomposable symmetric representations to be
equal to  $\Delta(\Gamma_{\pi,\sigma})$.

As for quivers and generalized quivers, we define finite and tame signed quivers and
classify them (\ref{finite_class},\ref{tame_class}).
We describe the indecomposable representations of finite and tame quivers. 
By \ref{embedding}, for any signed quiver this description can be reduced to 
the same for a quiver such that $\pi$ is admissible. 
Considering the finite and tame quivers case by case we
construct the indecomposable symmetric representations explicitly
and check that the condition of \ref{dim_sym_rep} is fulfilled.
Thus we get in \ref{tame_dim} that five tame signed
quivers yield the root systems corresponding to the graphs:
$D_n^{(2)},C_n^{(1)},A_{2n}^{(2)},B_n^{(1)},A_{2n-1}^{(2)}$.

\section{Signed quivers.}

\begin{definition}\label{signed_quiver}
A signed quiver is a triple $Q^{\sigma}=(Q,\sigma,m)$, where $Q$ is a quiver, 
$\sigma: Q_0\cup Q_1 \to \{-1;0;1\}$ is a sign function, 
$Q_0 = \{1,1^*,\cdots,m,m^*,m+1,\cdots,n\}$ 
contains $m\geq 0$ pairs of vertices that we call {\it twins}. 
We also call a vertex $i$ {\it signed} if $\sigma(i)\neq 0$; otherwise
we call it unsigned. The same terms apply to the arrows. The data is subject
to the axioms:

{\bf 1.} twins are unsigned

{\bf 2.} a non-loop arrow is signed if and only if its vertices are twins

{\bf 3.} a loop is signed if and only if its vertex is signed.
\end{definition}

To draw a signed quiver, we draw the underlying quiver and provide the signed
vertices and arrows with their signs; vertices connected by signed arrows 
are assumed to be twins. For a pair of twins that are not connected, we have
to index them explicitly as $i$ and $i^*$ for some $i$.

\begin{definition}
We say that a dimension $\alpha\in {\bf Z}_+^{Q_0}$ is sign-matched,
if $\alpha_i = \alpha_{i^*}$ for twins $i$ and $i^*$ and 
$\alpha_j$ is even, if $\sigma(j)=-1$.
\end{definition}

\begin{definition}\label{rep_sign}
Let $Q^{\sigma}$ be a signed quiver, let $\alpha$ be a sign-matched
dimension. Fix the data as follows:

$\bullet$ a vector space $V(i)$ for each $i\in Q_0$ such that for twins $i,i^*$, $V(i^*)=V(i)^*$      

$\bullet$ for any signed $k\in Q_0$, a linear map $J_k:V(k)\to V(k)^*$  such that  $J_k^*=\sigma(k)J_k$.

\noindent Denote by $G(Q^{\sigma},\alpha)$ the direct product of groups as follows: 

$GL(V(k))$, for every single unsigned vertex $k$ and every pair $k,k^*$ of twins,

$O(V(i))$, for every $i\in Q_0$ with $\sigma(i)=1$,

$Sp(V(j))$, for every $j\in Q_0$ with $\sigma(j)=-1$.

\noindent Denote by $R(Q^{\sigma},\alpha)$
          the direct sum of $G(Q^{\sigma},\alpha)$-modules as follows:

${\rm Hom}(V(i),V(j))$, for an unsigned arrow $i\to j$,

$S^2 V(i)^*$, for $i\xrightarrow{+} i^*$ or a $+$ loop on signed $i$; 
$S^2 V(i)$,   for $i\xleftarrow{+}  i^*$,

$\wedge^2 V(i)^*$, for $i\xrightarrow{-} i^*$ or a $-$ loop on signed $i$;  
$\wedge^2 V(i)$, for $i\xleftarrow{-}  i^*$.
\end{definition}

The group $(G(Q^{\sigma},\alpha),R(Q^{\sigma},\alpha))$ is in general neither connected
nor semi-simple. Denote by $(S(Q^{\sigma},\alpha),R(Q^{\sigma},\alpha))$ 
the commutant of its unity connected component (this is a result of replacing the factors
$GL(V(i))$ and $O(V(j))$ of $G(Q^{\sigma},\alpha)$ by $SL(V(i))$ and $SO(V(j))$, respectively). 
It is not difficult to 
show that the set of connected 
semi-simple groups that we obtain this way can be defined differently,
as follows:
\begin{proposition} Let $(H,U)$ be a connected semi-simple linear group such that
each irreducible $H$-factor of $U$ is isomorphic (as a linear group) to one
of the items:
\begin{equation}
\begin{split}
        (SL(L)\times SL(M),{\rm Hom}(L,M) ),
        & (SL(L)\times SO(M),{\rm Hom}(L,M) ),\\
        (SL(L)\times Sp(M),{\rm Hom}(L,M) ),
        & (SL(L),\wedge^2 L), (SL(L), S^2 L), \\
       (SO(L)\times SO(M),{\rm Hom}(L,M) ),
       & (SO(L)\times Sp(M),{\rm Hom}(L,M) ),\\
        (Sp(L)\times Sp(M),{\rm Hom}(L,M) ),
       & (SO(L),\wedge^2 L), (Sp(L), S^2 L) \\
        (SL(L),{\rm End}_0(L) ),
       & (Sp(L),\wedge^2_0 L), (SO(L), S^2_0 L). 
\end{split}
\end{equation}
Then there exist a signed quiver $Q^{\sigma}$ and a sign-matched dimension $\alpha$ such that
$(H,U/U^H)\cong (S(Q^{\sigma},\alpha),R(Q^{\sigma},\alpha)/R(Q^{\sigma},\alpha)^{S(Q^{\sigma},\alpha)})$.
\end{proposition}

\begin{remark}\label{DW-gq}
The signed quivers generalize the {\it generalized quivers}
from \cite{dw1}. Namely, the generalized quivers
are the groups $(G(Q^{\sigma},\alpha),R(Q^{\sigma},\alpha))$ for
the signed quivers of the special sort:

orthogonal case:  all signed arrows are - and all signed vertices are +;

symplectic case:  all signed arrows are + and all signed vertices are -.  
\end{remark}

\begin{remark}\label{nice_prop}
One of the natural tasks of Invariant Theory consists in classifying actions 
of reductive groups with nice properties of orbits and invariants.
Examining the known classifications, 
one can see that an important part of the lists consists of the groups 
arising from signed quivers, especially in what concerns "serial" cases.
For example the connected irreducuble linear groups with finitely many orbits are classified 
in \cite{kac}, Theorem 2; one can observe that all infinite series of groups arise from signed quivers 
except for the unique example: 
$$(SL_2\times SL_3\times SL_n, {\bf k}^2\otimes {\bf k}^3\otimes {\bf k}^n),
 n\geq 3.$$
Analogously, almost all serial irreducible semi-simple linear 
groups with polynomial algebra of invariants (see \cite{lit}) arise from signed quivers.
\end{remark}


\section{Signed quivers and symmetric quivers.}

The subsequent definition follows that of Derksen and Weyman:

\begin{definition}\label{sym_quiv}
(cf. \cite{dw1}) A symmetric quiver $S$ is a quiver endowed with
an involution $*$ acting on both $S_0$ and $S_1$ such that
$$t(\varphi^*) = (h\varphi)^*, h(\varphi^*) = (t\varphi)^*.$$ 
\end{definition}

\noindent For a signed quiver $Q^{\sigma}$, we first define 
$*$ on subsets of $Q_0$ and $Q_1$:
$i^*= i^*, (i^*)^* = i$ for a pair $i,i^*$ of twins, $k^*=k$
for a signed vertex $k$; $\varphi^*=\varphi$, for a signed arrow
$\varphi$. To define $*$ on the whole of $Q$ we
have to extend $Q$, as follows: 
we add a vertex $j^*$ for any 
single unsigned vertex $j$ of $Q$ and an
arrow $\varphi^*:q^*\to p^*$
for any unsigned arrow $\varphi:p\to q$ of $Q^{\sigma}$. 
We denote by $\widetilde{Q}$ the quiver that we 
obtained this way. Clearly, the defined involution $*$ 
fulfills definition \ref{sym_quiv}.
\begin{example}
Let $Q^{\sigma}$ be as follows:
${\circ} \xrightarrow{} {\circ} \xleftarrow{} {\circ} \xrightarrow{-} {\circ}$.
We numerate the vertices from left to right as $1,2,3,3^*$. Then $\widetilde{Q}$
has vertices $1,2,3,3^*,2^*,1^*$ and looks like: 
${\circ} \xrightarrow{} {\circ} \xleftarrow{} {\circ} \xrightarrow{-} {\circ}
 \xleftarrow{} {\circ} \xrightarrow{} {\circ}$.
\end{example}
We therefore constructed a symmetric quiver $\widetilde{Q}$, which is also signed 
such that all signed vertices and arrows belong to $Q^{\sigma}$. 
The key tool in dealing with the representations of $\widetilde{Q}$ is the duality functor.
For $\varphi\in \widetilde{Q}_1$ set $s_{\varphi} = -1$ 
if either $\sigma(\varphi)=-1$ and $\sigma(h\varphi)\neq -1$ 
or $\sigma(\varphi)\neq -1$ and $\sigma(h\varphi)=-1$,
and $s_{\varphi} = 1$, otherwise.
\begin{definition}\label{v^*}
Let $V$ be a representation of $\widetilde{Q}$.
The dual representation $V^*$ is:
$$V^*(i) = (V(i^*))^*, i\in S_0,$$
$$V^*(\varphi) = s_{\varphi} V(\varphi^*)^*: 
(V(h\varphi^*))^*\to V(t\varphi^*)^*, 
\varphi\in S_1.$$ 
\end{definition}
\begin{remark}
The sign $*$ is used in the above definition
in three different senses: for the involution of $\widetilde{Q}$, 
for the dual vector space, and for the dual linear map.
\end{remark}
\begin{remark}
Our definition of the dual representation differs from
that from \cite{dw1} by the multiple $s_{\varphi}$.
This difference allows us to recover a wider set of groups.
\end{remark}
\begin{remark}\label{v**}
 Clearly, we have: $(V^*)^*(\varphi) = s_{\varphi}s_{\varphi^*}V(\varphi)$,
and $s_{\varphi}s_{\varphi^*}=-1$ when exactly one of $t\varphi,h\varphi$ is a $-$ vertex.
So $(V^*)^*$ is not always the same as $V$. However  $(V^*)^*$ is always conjugate to $V$ by
the product of $-{\rm Id}$ operators taken over all $-$ vertices of  $Q^{\sigma}$.
\end{remark}

\begin{definition}\label{si-mi}
We say that a dimension $\gamma\in {\bf Z}_+^{\widetilde{Q}_0}$ is sign-matched,
if $\gamma_i = \gamma_{i^*}$ and $\gamma_j$ is even, if $\sigma(j)=-1$.
\end{definition}

Clearly, extending by symmetricity a sign-matched dimension $\alpha$ for $Q$ 
to a dimension vector $\tilde{\alpha}\in {\bf Z}_+^{\widetilde{Q}_0}$
we get $\tilde{\alpha}$ is sign-matched.
We define an embedding of the space $R(Q^{\sigma},\alpha)$
to $R(\widetilde{Q},\tilde{\alpha})$.
Recall that the definition of $R(Q^{\sigma},\alpha)$ is
based, in particular on fixing vector spaces $V(i),i\in Q_0$
such that $V(i^*)=V(i)^*$ for a pair $i,i^*$ of twins,
and on a choice of the isomorphism
$J_k:V(k)\to V(k)^*$ for any signed vertex $k$.
For any vertex $i\in\widetilde{Q}_0$ set
$V(i)={\bf k}^{\tilde{\alpha}_i}$ and choose a structure as follows:
\begin{definition}\label{sign_form}
A {\it signed form} on the spaces $V(i),i\in \widetilde{Q}_0$
is a collection of isomorphisms
$J_i: V(i)\to (V(i^*))^*, i\in\widetilde{Q}_0$,
such that $J_{i^*}= (J_i)^*$, if $i\neq i^*$, and
$(J_i)^*=\sigma(i)J_i$, otherwise.
\end{definition}
\begin{proposition}\label{allforms}
 All signed forms are conjugate by $GL(\tilde{\alpha})$.
\end{proposition}
\begin{proof}
This follows from the fact that all non-degenerate symmetric and  anti-symmetric 
forms on a vector space $W$ are conjugate by $GL(W)$ and (since $J_{i^*}$ is determined
by $J_i$) that all non-degenerate maps $U\to W^*,\dim U = \dim W$ are conjugate by 
$GL(U)\times GL(W)$.
\end{proof}

So all signed forms are equivalent to each other and once a form is chosen, we are able 
to introduce our main concept:
 
\begin{definition}\label{sym_rep}
A representation $V\in R(\widetilde{Q},\tilde{\alpha})$
is called symmetric, if the collection of maps
$\{J_i:V(i)\to V^*(i),i\in \widetilde{Q}_0\}$ is
an isomorphism of $V$ to $V^*$.
\end{definition}

\begin{example}\label{+-}
Let $Q^{\sigma}$ (and also $\widetilde{Q}$) consist of one + vertex $a$ and one - loop $\varphi$
on it. By Definition \ref{v**}, we have for a representation $V$: $V^*(\varphi) = - V(\varphi)^*$.
A signed form
$J_a$ is in this case just a symmetric form on $V(a)$. The symmetric representations 
are subject to the equation: $V(\varphi)^*J_a+J_aV(\varphi)=0$ or for $A= J_aV(\varphi): A^* = -A$,
hence, these are the anti-self-adjoint endomorphisms of $V(a)$.
Note that the set of symmetric representations is a $O(V(a))$-stable vector subspace in
${\rm End}(V(a))$ isomorphic to $\wedge^2 V(a)$. 
\end{example}
\begin{example}\label{++}
Let $Q^{\sigma}$ (and also $\widetilde{Q}$) consist of one + vertex $a$ and one + loop $\varphi$
on it. The symmetric representations are the self-adjoint endomorphisms and constitute
a $O(V(a))$-module $S^2 V(a)$.
\end{example} 
\begin{example}\label{wedge}
Let $Q^{\sigma}$ (and also $\widetilde{Q}$) be ${\circ} \xrightarrow{-} {\circ}$.
Denote by $\varphi$ the unique arrow, $a=t\varphi,b=h\varphi$. 
By Definition \ref{v**}, we have for a representation $V$: $V^*(\varphi) = - V(\varphi)^*$. 
The symmetric representations 
are subject to the equation: $V(\varphi)^*J_a+J_bV(\varphi)=0$ or for $A= J_bV(\varphi): A^* = -A$.
So the set of symmetric representations is naturally isomorphic to $\wedge^2 V(a)^*$.
\end{example}

We define an involution $\tau\in GL(R(\widetilde{Q},\tilde{\alpha}))$;
namely, for any arrow $\varphi\in \widetilde{Q}_1$ set 
\begin{equation} 
\tau(V)(\varphi) = 
J_{h\varphi}^{-1} V^*(\varphi) J_{t\varphi}\in 
{\rm Hom}(V(t\varphi),V(h\varphi)).
\end{equation}
\noindent Clearly, a representation $V\in R(\widetilde{Q},\tilde{\alpha})$
is symmetric if and only if $\tau(V) = V$.
\begin{proposition} $\tau ^2 = {\rm Id}$.
\end{proposition}
\begin{proof}
$\tau^2(V)(\varphi) = J_{h\varphi}^{-1} ((\tau(V))^*(\varphi) J_{t\varphi} = 
J_{h\varphi}^{-1} s_{\varphi} ( \tau(V)(\varphi^*) )^* J_{t\varphi} = $
$$J_{h\varphi}^{-1} s_{\varphi} (J_{h\varphi^*}^{-1} 
s_{\varphi^*}(V(\varphi))^*J_{t\varphi^*} ) )^* J_{t\varphi} = 
s_{\varphi}s_{\varphi^*} J_{h\varphi}^{-1} J_{(h\varphi)^*}^* V(\varphi)
J_{(t\varphi)^*}^{-*} J_{t\varphi},$$
\noindent (here and below we write $X^{-*}$ instead of $(X^*)^{-1}=(X^{-1})^*$).
Note that $J_{h\varphi}^{-1} J_{(h\varphi)^*}^*$ (resp. $J_{(t\varphi)^*}^{-*} J_{t\varphi}$)
is the $\pm 1$ scalar operator, -1 iff  $\sigma(h\varphi)=-1$ (resp. $\sigma(t\varphi)=-1$).
By Remark \ref{v**}, $s_{\varphi}s_{\varphi^*}=-1$ when exactly one of $t\varphi$
and $h\varphi$ is a $-$ vertex.
\end{proof}

\begin{proposition}\label{inner}
For $g\in GL(\tilde{\alpha})$, $g^{\tau} = \tau g \tau$ belongs to the image of 
$GL(\tilde{\alpha})$ in $GL(R(\widetilde{Q},\tilde{\alpha}))$,
$g^{\tau}_i =   
J_i^{-1} g_{i^*}^{-*} J_i$.
\end{proposition}
\begin{proof}
First consider an unsigned arrow $\varphi\in Q_1, t\varphi = k, h\varphi = l$.
Set $A=V(\varphi),B=V(\varphi^*)$. Then we have:
$$
\tau(A,B)=(s_{\varphi}J_l^{-1} B^* J_k, s_{\varphi^*}J_{k^*}^{-1} A^* J_{l^*}),
g(A,B) = (g_l A g_k^{-1}, g_{k^*} B g_{l^*}^{-1}).
$$
$$\tau g \tau(A,B) = (
s_{\varphi} s_{\varphi^*} J_l^{-1} g_{l^*}^{-*} J_{l^*}^* A J_{k^*}^{-*} g_{k^*}^* J_k,
s_{\varphi^*} s_{\varphi} J_{k^*}^{-1} g_k^{-*} J_k^* B J_l^{-*} g_l^* J_{l^*}),$$
and we are done by Remark \ref{v**}.
Similarly, we consider a signed $\varphi$. 
\end{proof}

\begin{corollary}\label{centralizer}
The centralizer $GL(\tilde{\alpha})^{\tau}$ of $\tau$ is generated by
the kernel for the action of $GL(\tilde{\alpha})$ and the subgroup
\begin{equation}\label{Gtau}
G^{\tau}=\{ g\in GL(\tilde{\alpha})\vert J_i^{-1} g_{i^*}^{-*} J_i = g_i,
i\in \widetilde{Q}_0\}.
\end{equation}
\end{corollary}

\begin{corollary}\label{iso_of_group}
The linear groups $(GL(\tilde{\alpha})^{\tau},
R(\widetilde{Q},\tilde{\alpha})^{\tau} )$ and
$(G(Q^{\sigma},\alpha),R(Q^{\sigma},\alpha))$
are isomorphic.  
\end{corollary}
\begin{proof}
By Corollary \ref{centralizer}, the groups $(GL(\tilde{\alpha})^{\tau},
R(\widetilde{Q},\tilde{\alpha})^{\tau} )$ and 
$(G^{\tau},R(\widetilde{Q},\tilde{\alpha})^{\tau} )$ are equal.
Clearly, $G^{\tau}$ is naturally isomorphic to $G(Q^{\sigma},\alpha)$, so we only
need to compare Definition \ref{rep_sign} of $R(Q^{\sigma},\alpha)$ and
Definition \ref{sym_rep}. By Definition \ref{sym_rep}, the map $V(\varphi^*)$ is uniquely defined
by $V(\varphi)$ for any symmetric representation $V$ and any unsigned
arrow $\varphi\in Q_1$. Since this is the unique restriction for $V(\varphi)$ and $V(\varphi^*)$,
we get the submodule ${\rm Hom}(V(t\varphi),V(h\varphi))\subseteq R(\widetilde{Q},\tilde{\alpha})^{\tau} $.
If $\varphi$ is signed, then in Examples \ref{+-}, \ref{++}, \ref{wedge}
we saw that the symmetricity condition of $V(\varphi)$ gives rise to a submodule as in 
Definition \ref{rep_sign}.
\end{proof}

\begin{theorem}\label{isofiso} Symmetric representations are conjugate by $G^{\tau}$ if and
only if they are isomorphic.
\end{theorem}
\begin{proof}
The proof is based on an observation from \cite{mwz}, as follows.
Let $X$ be a set acted upon by a group $G$ and an involution
$\sigma$, which normalizes $G$. Assume furthermore that
$G$ is a subgroup in the group of invertible elements of a
finite-dimensional associative algebra $A$ 
and the anti-automorphism $g\to \sigma g^{-1} \sigma$ of $G$ extends
to a ${\bf k}$-linear involution of $A$. 

\begin{proposition}\label{zele}
\cite{mwz}  Suppose that for any $x\in X$ each invertible element of the linear span
in $A$ of the stabilizer $G_x$ belongs to $G_x$.
Then any two points in $X^{\sigma}$ are $G$-conjugate
if and only if they are $Z_G(\sigma)$-conjugate.
\end{proposition} 

We apply Proposition \ref{zele} to $X=R(\widetilde{Q},\tilde{\alpha})$,
$A=L(\tilde{\alpha})$, $G=GL(\tilde{\alpha})$ and $\sigma=\tau$.
Actually, by Proposition \ref{inner}, $(\tau g^{-1} \tau)_i =  J_i^{-1} g_{i*}^* J_i$
hence, this anti-automorphism of $GL(\tilde{\alpha})$ extends to the whole of 
$L(\tilde{\alpha})$. For any $V\in X$, $G_V={\rm Aut}(V)$ 
and its linear span is ${\rm End}(V)$. It remains to note that by Corollary
\ref{centralizer}, the action of $G^{\tau}$ on $R(\widetilde{Q},\tilde{\alpha})$
is equal to that of $Z_G(\sigma)=GL(\tilde{\alpha})^{\tau}$. 
\end{proof}

\begin{definition}
Let $U,V$ be symmetric representations of $\widetilde{Q}$ 
with respect to signed form $J^U,J^V$, respectively.
We say that $U$ and $V$ are {\it symmetrically isomorphic}
if there exists an isomorphism $T:U\to V$ such that
$T$ is an isometry of the signed forms: for any $x\in U(i), y\in U(i^*)$ holds:
$\langle J_i^U(x),y\rangle =  \langle J_i^V(Tx),Ty\rangle$.
\end{definition}
Proposition \ref{allforms} and Theorem \ref{isofiso} imply:
\begin{corollary}\label{outeriso}
Symmetric representations are symmetrically isomorphic if and only 
if they are isomorphic.
\end{corollary}

\section{Indecomposable symmetric representations.}

In this section we introduce and describe  the indecomposable symmetric
representations of the symmetric quiver $\widetilde{Q}$ corresponding
to a signed quiver $Q^{\sigma}$. Again these results 
and their proofs generalize those from \cite{dw1}.

Let $V^i$ be a symmetric representation of $\widetilde{Q}$ 
with respect to a signed form $J^i$, $i=1,2$. Then $V_1+V_2$ is symmetric
with respect to the signed form $J=J^1\oplus J^2$.  

\begin{definition}\label{indecomp-sym}
A symmetric representation $V$ of $\widetilde{Q}$ is called indecomposable if
$V$ is not isomorphic to a non-trivial direct sum $W=U+V$ of symmetric representations.
\end{definition}

By Corollary \ref{outeriso}, a symmetric representation is determined, modulo 
symmetric isomorphism by the isomorphism classes of indecomposable
symmetric summands of any splitting. 
It is however not clear  that all splittings of a representation
into indecomposable ones have the same summands modulo permutations and 
isomorphisms, as it is for the usual quivers.
It follows directly from Definition \ref{sym_rep} 
of symmetric representations:

\begin{def-prop}\label{orthog}
Let $V$ be a symmetric representation of $\widetilde{Q}$,
$W\subseteq V$ be a subrepresentation. Set 
$W^{\perp}_i = J_i^{-1}(W_{i^*}^{\perp})\subseteq V_i$,
where $W_{i^*}^{\perp}\subseteq V_{i^*}^*$ is the 
annulator.
This collection of subspaces is a subrepresentation in $V$.
\end{def-prop}

\begin{proposition}\label{sym-split}
Let $V$ be a symmetric representation of $\widetilde{Q}$,
$W\subseteq V$ be a subrepresentation. If 
$W_i\cap W^{\perp}_i = 0$ for any 
$i\in\widetilde{Q}_0$, then both $W$ and $W^{\perp}$
carry structures of symmetric representations such
that $V$ splits: $V=W+W^{\perp}$.
\end{proposition}

\begin{proof}
Assuming $W_i\cap W^{\perp}_i = 0$, we
get an inequality $\dim W_i +\dim W^{\perp}_i \leq \dim V_i$, which 
can be rewritten as $\dim W_i \leq \dim W_{i^*}$ (because $\dim V_i= \dim V_{i^*}$).
These inequalities for $i$ and $i^*$ yield: $\dim W_i = \dim W_{i^*}$.
Hence, $\dim W_i + \dim W^{\perp}_i = \dim V_i$ and  $V=W+W^{\perp}$ splits 
as a representation of $\widetilde{Q}$.
Furthermore, the condition $W_i\cap W^{\perp}_i = 0$ implies
$J_i W_i$ is naturally isomorphic to $W_{i^*}^*$. Hence, restricting
the form $J_i$ to $W_i$ for all $i$, we get
a signed form on $W$. Similarly, we endow $W^{\perp}$ with a structure 
of symmetric representation.
\end{proof}
There is a natural way to construct symmetric representations, as follows:
\begin{proposition}\label{v+v*}
For any representation $V$ of $\widetilde{Q}$, $V+V^*$ carries a
natural structure of symmetric representation.
\end{proposition}
\begin{proof}
For $X=V+V^*$, $X(i) = V(i)\oplus V(i^*)^*$. Hence, $(X(i^*))^* = X(i)$
and we set $J_i = {\rm Id}$ for $i\neq i^*$; clearly $J_i^* = J_{i^*}$. 
If $i$ is signed, then we define 
$$J_i:V(i)\oplus V(i)^*\to V(i)^* \oplus V(i), J_i(v+v^*) =\sigma(i)v^* + v$$
for any $v\in V(i),v^*\in V(i)^*$. 
In this case we have: $J_i^*(v+v^*)=v^*+\sigma(i)v=\sigma(i) J_i$.
A straightforward calculation shows that  $J$ is an isomorphism of $X$ to $X^*$. 
\end{proof}

\begin{lemma}\label{1_or_2}
Let $V$ be an indecomposable symmetric representation of $\widetilde{Q}$.
Then as a representation of $\widetilde{Q}$
$V$ is either indecomposable or $V=W+W^*$ for an indecomposable $W$. 
\end{lemma}
\begin{proof}
Assume that $V$ splits as a representation of $\widetilde{Q}$: $V=W+U$ such that $W$ 
is indecomposable. By Proposition \ref{sym-split}, $W_i\cap W^{\perp}_i \neq 0$ 
for some $i$. 

Note that we have another splitting:
$V= W^{\perp} +U^{\perp}$.
First we claim: $W_i\cap U^{\perp}_i = 0$ for any $i$.
Consider two maps:
\begin{equation}
P_1:U^{\perp}\hookrightarrow V = W +U 
\stackrel{p_1}{\longrightarrow} W, \quad
P_2:W\hookrightarrow V = W^{\perp} +U^{\perp} 
\stackrel{p_2}{\longrightarrow} U^{\perp} 
\end{equation}
Clearly, $P_1$ and $P_2$ are homomorphisms of representations
of $\widetilde{Q}$, hence, $P_1P_2$ is an endomorphism of $W$.
The kernel of $P_2$ is equal to $W\cap W^{\perp}$,
hence, is nonempty. So $P_1P_2$ is not invertible, hence,
should be nilpotent, because the endomorphism ring
of the indecomposable representation $W$ is local.
On the other hand, any element of $W_i\cap U^{\perp}_i$
is stable under both $P_1$ and $P_2$. Thus
$W_i\cap U^{\perp}_i = 0$ for any $i$.

Note that, as representation of $\widetilde{Q}$,
$U^{\perp}\cong (V/U)^*\cong W^*$. Then the dimension of
$X=W+U^{\perp}$ is symmetric. We claim: $X_i\cap X^{\perp}_i = 0$
for all $i$.
Consider a map $T: X\to X$ given by the matrix
$$\left(
\begin{array}{cc}
{\rm Id_W} & P_1 \\
P_2 & {\rm Id}_{U^{\perp}}
\end{array}
\right)$$
Taking into account the equality
$X^{\perp}_i = W_{i}^{\perp} \cap U_i$, one can easily 
check that the kernel of $T$ is exactly $X\cap X^{\perp}$.
That $P_1P_2$ is nilpotent implies the same for $T-{\rm Id}_X$,
hence, $T$ is invertible and the claim is proved.
By Proposition \ref{sym-split}, $V$ splits as $V=X+Y$,
so $V=X=W+U^{\perp}\cong W+W^*$.
\end{proof}

\begin{proposition}\label{2_not_1}
If $W$ is an indecomposable representation of $\widetilde{Q}$
such that there exists a symmetric representation $W'\cong W$,
then the symmetric representation $W+W^*$ is decomposable.
\end{proposition}
\begin{proof}
By assumption, $\dim W=\dim W^*$. Hence, we have two symmetric
representations of the same dimension: $2W'$ and $W+W^*$.
Since these are isomorphic as representations of $\widetilde{Q}$,
by Theorem \ref{isofiso}, they are symmetrically isomorphic. 
Therefore, $W+W^*$ splits.
\end{proof}

\begin{theorem}\label{class_indec}
For any symmetric representation $V$ of $\widetilde{Q}$ there exists a unique,
modulo isomorphisms and permutations of factors 
decomposition of $V$ into a direct sum of indecomposable symmetric
representations.
\end{theorem}
\begin{proof}
By Lemma \ref{1_or_2},
each splitting of $V$ into a sum of indecomposable symmetric
representations consists of items $W_1$ or $W_2+W_2^*$, where
$W_1$ and $W_2$ are indecomposable representations of $\widetilde{Q}$.
By Proposition \ref{2_not_1}, none $W_2$-like summand can be isomorphic
to a $W_1$-like one. Hence, the summands of all splittings
are isomorphic as representations of $\widetilde{Q}$, modulo permutations. 
By Theorem \ref{isofiso}, they are also symmetrically isomorphic.
\end{proof}
\section{Two functors.}
\begin{definition}\label{del}
Let $Q^{\sigma}$ be a signed quiver and $i\in Q_0$ a signed
vertex incident to a unique arrow $\varphi$; denote by $j$ the
other vertex of $\varphi$ and assume $j\neq i$. Replacing $i$ by $j^*$ and
$\varphi$ by a $\sigma(i)$-signed arrow $\psi$ between $j$ and $j^*$ we
get a new signed quiver that we denote by ${\rm Del}_i(Q)$.
For a representation $V$ of $\widetilde{Q}$ denote by ${\rm Del}_i(V)$ 
the representation of $\widetilde{{\rm Del}_i(Q)}$ such that
${\rm Del}_i(V)(\psi)$ is the composition (in the appropriate order)
of $V(\varphi)$ and $V(\varphi^*)$. 
\end{definition}

\begin{example}
$Q^{\sigma}:{\circ}_j\rightarrow\overset{+}{\circ}_i$\quad
$\widetilde{Q}:{\circ}_j\rightarrow\overset{+}{\circ}_i\rightarrow{\circ}_{j^*}$\quad
${\rm Del}_i(Q)=\widetilde{{\rm Del}_i(Q)}:{\circ}_j\xrightarrow{+}{\circ}_{j^*}$.
\end{example}

Since the vertices of $\widetilde{{\rm Del}_i(Q)}$ are vertices of $\widetilde{Q}$,
a signed form for $\widetilde{{\rm Del}_i(Q)}$ can be inherited from $\widetilde{Q}$.
One can easily deduce from definitions:
\begin{proposition} If $V$ is a symmetric representation, then ${\rm Del}_i(V)$ is.
\end{proposition}

\begin{definition}\label{ins}
Let $Q^{\sigma}$ be a signed quiver and $\varphi\in Q_1$ a signed
non-loop arrow. Inserting a new $\sigma(\varphi)$-signed vertex $a$
and replacing the arrow $\varphi$ by an unsigned arrow $\psi:t\varphi\to a$ we get
a new signed quiver that we denote by ${\rm Ins}_{\varphi}(Q)$.
For a representation $V$ of $\widetilde{Q}$ denote by ${\rm Ins}_{\varphi}(V)$ 
the representation of $\widetilde{{\rm Ins}_{\varphi}(Q)}$ such that
${\rm Ins}_{\varphi}(V)(a)={\rm Im}V(\varphi)$, ${\rm Ins}_{\varphi}(V)(\psi)=V(\varphi)$,
and ${\rm Ins}_{\varphi}(V)(\psi^*)$ is the embedding ${\rm Im}V(\varphi)\to V(h\varphi)$.  
\end{definition}

\begin{example}
$Q^{\sigma}=\widetilde{Q}:{\circ}\xrightarrow{\varphi,+}{\circ}$\quad
${\rm Ins}_{\varphi}(Q):{\circ}\xrightarrow{\psi}\overset{+}{\circ}_a$\quad
$\widetilde{{\rm Ins}_{\varphi}(Q)}:{\circ}\xrightarrow{\psi}
\overset{+}{\circ}_a\xrightarrow{\psi^*}{\circ}$.
\end{example}

\begin{proposition}
Set $p={\rm Ins}_{\varphi}(V)(\psi)$, $i={\rm Ins}_{\varphi}(V)(\psi^*)$.
If $V$ is a symmetric representation, then the  map $J_a=i^* J_l p^{-1}:V(a)\to V(a)^*$ 
is well defined, $J_a$ and the signed form $J$ of $\widetilde{Q}$ yield a signed form
for $\widetilde{{\rm Ins}_{\varphi}(Q)}$, and ${\rm Ins}_{\varphi}(V)$ is symmetric.
\end{proposition}

\begin{proof} Set $l=t\varphi,l^*=h\varphi$, $A=V(\varphi)=ip$.
Consider a diagram as follows:
$$
\begin{CD}
V(l)         @>p>>   V(a)    @>i>>  V(l^*) \\
@VJ_lVV               @VVi^* J_l p^{-1}V    @VV{J_{l^*}}V\\
V(l^*)^*     @>>i^*> V(a)^*  @>>\sigma(\varphi)p^*> V(l)^*
\end{CD}
$$
Since $p$ is surjective by definition, $p^{-1}$ is defined
on the whole of $V(a)$, uniquely modulo the kernel of $p$. 
Assume that $V$ is symmetric: $\sigma(\varphi)A^* J_l=J_{l^*}A$.
If $p(v)=0$, then $A(v)=0$ and 
$\sigma(\varphi)p^* i^* J_l(v)=\sigma(\varphi)A^* J_l(v)=J_{l^*}A(v)=0.$
Since $p^*$ is injective, this implies $i^* J_l(v)=0$, so $J_a=i^* J_l p^{-1}$ 
is well defined. It is also clear that both squares are commutative.
Finally, $J_a^*=p^{-*}J_l^*i$, and since the right square is commutative,
$J_l^*i=\sigma(\varphi)p^*J_a$, so $J_a^*=\sigma(\varphi)J_a$.
\end{proof}

Remark that the functors ${\rm Del}$ and ${\rm Ins}$ are inverse to each other
for quivers: ${\rm Del}_a({\rm Ins}_{\varphi}(Q)) = Q$, ${\rm Ins}_{\psi}({\rm Del}_i(Q)) = Q$.
One can easily prove:

\begin{proposition}\label{vice-versa}
{\bf 1}. Assume that the quiver $Q$ is as in {\rm \ref{ins}} and $V$ is a symmetric representation of $\widetilde{Q}$.
Then ${\rm Del}_a({\rm Ins}_{\varphi}(V)) = V$. 

{\bf 2}. Assume that the quiver $Q$ is as in {\rm \ref{del}}; set $\gamma = \varphi$ or $\gamma = \varphi^*$
such that $h\gamma = i$. If $V$ is a symmetric representation 
of $\widetilde{Q}$ such that $V(\gamma)$ is surjecive, then ${\rm Ins}_{\psi}({\rm Del}_i(V)) = V$.
\end{proposition}

\begin{corollary}
Assume that the quiver $Q$ is as in {\rm \ref{ins}}. Then the functor ${\rm Ins}$ is an
isomorphism of the symmetric representations of $\widetilde{Q}$ onto those symmetric representations $W$
of $\widetilde{{\rm Ins}_{\varphi}(Q)}$ such that $W(\psi)$ is surjective.
\end{corollary}

\begin{proposition}\label{embedding}
Assume that the quiver $Q$ is as in \ref{ins}. Then the 
indecomposable symmetric representations of $\widetilde{Q}$ 
are the representations ${\rm Del}_a(W)$ for all indecomposable symmetric
representations $W$ of $\widetilde{{\rm Ins}_{\varphi}(Q)}$ such that $W(\psi)$ is surjective.
\end{proposition}
\begin{proof}
Clearly both functors Del and Ins respect direct sums of representations. Hence,
by Proposition \ref{vice-versa}.2, ${\rm Del}_a(W)$ is indecomposable for an
indecomposable representation $W$,
and by Proposition \ref{vice-versa}.1, all indecomposable symmetric representations of 
$\widetilde{Q}$ can be obtained this way.
\end{proof}
\section{Root systems for symmetrizable generalized Cartan matrices.}
The notion of a generalized Cartan matrix (GCM) and the associated Weyl group and
root system was introduced by V.Kac (\cite{kac2}, \cite{kac3}).
In \cite{kac2} is proved that the root system corresponding
to a {\it symmetric} GCM $A$ is exactly the dimensions of the indecomposable
representations of a quiver having the Dynkin diagram $S(A)$ of $A$ as the underlying graph.
Moreover, $\alpha$ is a real root if and only if there exists a unique indecomposable
of dimension $\alpha$, modulo isomorphism. 

It is of interest to consider not only symmetric GCM but also the {\it symmetrizable} ones. 
In particular, all GCM of finite and affine types are symmetrizable (see their classification in 
\cite{kac2}, \cite{kac3}).
In some cases the root system corresponding to a symmetrizable GCM can be
obtained factorizing (in a sense) the root system of a symmetric GCM $A$ by a diagram 
automorphism  of $S(A)$. This was used, in particular, in \cite[$\S$7.9]{kac3}
in order to construct the root systems for $B_l,C_l,F_4,G_2$ in terms of
those for $D_{l+1},A_{2l-1},E_6,D_4$, respectively. 
Recently A.Hubery observed in \cite{hu} that a description of this kind is possible
in the context of an {\it admissible} automorphism of the graph $S(A)$.
In particular, for a symmetric quiver $\widetilde{Q}$ such that 
the automorphism $*$ of the underlying graph is admissible (see \ref{from_sig})
this result gives a root system structure to 
the set consisting of the symmetric roots and the sums of different roots symmetric
to each other.
By Lemma \ref{1_or_2}, this set should be related
to that of dimensions of the indecomposable symmetric representations of
$\widetilde{Q}$, but only in case when all signed vertices are $+$.
In this section we suggest a more general approach such that 
for arbitrary signs we introduce a graph and a kind of symmetrization of roots,
both depending on signs, and prove that this symmetrization is the root system
corresponding to the graph (Lemma \ref{root_system}). Furthermore, we show in 
Corollary \ref{dim_sym_rep} that under some extra conditions this 
root system is precisely the set of dimensions of the indecomposable symmetric 
representations.

Let $\Gamma$ be a symply-laced graph and let $\pi$ be an involutive automorphism of $\Gamma$.
\begin{definition}\label{admis}
We say that $\pi$ is admissible if none edge of $\Gamma$ is $\pi$-stable.
\end{definition}
Assume from now on that $\pi$ is admissible and fix also a sign $\sigma(i)=\pm 1$ for any 
$i\in\Gamma^{\pi}$, the set of $\pi$-stable vertices of $\Gamma$. 

\begin{example}\label{from_sig}
Let $Q^{\sigma}$ be a signed quiver that has neither signed arrows nor arrows joining
signed vertices. Then the underlying graph $\Gamma$ of $\widetilde{Q}$ together with $\pi=*$
and the sign $\sigma$ yield a data as above.
\end{example}
Denote by $\Gamma_{\pi,\sigma}$ a graph with vertices corresponding to the  $\pi$-orbits
of the vertices of $\Gamma$ and the edges as follows (denote by $\overline{i}$ the $\pi$-orbit of $i$):

for an edge of $\Gamma$ between $i,j\notin \Gamma^{\pi}$ we draw a simple edge between 
$\overline{i}$ and $\overline{j}$ 

for an edge of $\Gamma$ between $i\notin \Gamma^{\pi}$ and $j\in \Gamma^{\pi}$ 
we draw a double edge between $\overline{i}$ and $\overline{j}$, oriented
to $j$ (resp. from $j$) if $\sigma(j)=1$ (resp. $\sigma(j)=-1$). 

\begin{example}
Let $Q^{\sigma}$ be $\overset{+}{\circ}\leftarrow{\circ}\rightarrow\overset{-}{\circ}$. Then $Q^{\sigma}$ meets
the condition of \ref{from_sig} and following the above definition we get a graph
${\circ}\Leftarrow{\circ}\Leftarrow{\circ}$.
\end{example}

We are now going to describe the root system corresponding to the graph $\Gamma_{\pi,\sigma}$
in terms of the root system $\Delta(\Gamma)$. Let $Q(\Gamma) = {\bf Z}^{\Gamma_0}$ be
the root lattice of $\Gamma$ and $Q(\Gamma)_+={\bf Z}_+^{\Gamma_0}$ be the positive vectors.
Consider a sublattice as follows:
$$Q(\Gamma)_{\pi,\sigma}=\{\alpha\in Q(\Gamma)\vert \alpha_{\pi(i)} = \alpha_i,\sigma(j)^{\alpha_j}=1,
j\in \Gamma^{\pi}\}.$$

\noindent In other words, we consider
the symmetric vectors having even coefficients on $-$ vertices. Note that by Definition \ref{si-mi},
in the situation of Example \ref{from_sig}, $Q(\Gamma)_{\pi,\sigma}$ is exactly
the set of sign-matched dimension vectors. Denote by $\Pi = \{\varepsilon_i,i=1,\cdots,n\}$ the
standard basis of $Q(\Gamma)$ (the simple roots).
For $i\notin\Gamma^{\pi}$ (resp. $\sigma(i)=1$, $\sigma(i)=-1$) set
$\beta_{\overline{i}}=\varepsilon_i+\varepsilon_{\pi(i)}$  
(resp. $\varepsilon_i$, $2\varepsilon_i$);
denote by $\Pi_{\pi,\sigma}$ the obtained set.
Clearly, $\Pi_{\pi,\sigma}$ is a basis of 
$Q(\Gamma)_{\pi,\sigma}$ and also a system of generators for the semi-group 
$Q(\Gamma)_{\pi,\sigma,+}=Q(\Gamma)_{\pi,\sigma}\cap Q(\Gamma)_+$.

Recall that the Tits form endows the lattice $Q(\Gamma)$ with a bilinear symmetric form 
$\langle \; , \, \rangle$ such that its matrix in the basis $\Pi$ is exactly 
the GCM matrix of $\Gamma$ divided by 2.
The group $W(\Gamma)$ generated by the simple reflections 
$r_i(\alpha) = \alpha - 2\langle\alpha ,\varepsilon_i\rangle \varepsilon_i$
is called the Weyl group. The elements of
$\Delta_+^{re}=W(\Gamma)\Pi\cap  Q(\Gamma)_+$ are called {\it real positive roots}.
Furthermore, the set $F_{\Gamma}$ of all $\alpha\in Q(\Gamma)_+$ such that 
$\langle\alpha ,\varepsilon_i\rangle \leq 0$ for all $i$ and the support of $\alpha$ is connected
is called the {\it fundamental domain}. The elements of 
$\Delta_+^{im}=W(\Gamma)F_{\Gamma}\subseteq Q(\Gamma)_+$
are called {\it imaginary positive roots}. The positive root system $\Delta(\Gamma)_+$
is by definition the union of  $\Delta_+^{re}$ and $\Delta_+^{im}$.

Now we construct the root system of $\Gamma_{\pi,\sigma}$ in $Q(\Gamma)_{\pi,\sigma}$.
For a vertex $\overline{i}$ of $\Gamma_{\pi,\sigma}$, let $s_{\overline{i}}\in W(\Gamma)$
be $r_i r_{\pi(i)}$ if $i\notin\Gamma^{\pi}$ and $r_i$, otherwise. Note that
 $r_i r_{\pi(i)}=r_{\pi(i)} r_i$, because $\pi$ is admissible. Let
$W(\Gamma)_{\pi}$ be the subgroup in $W(\Gamma)$ generated by $s_{\overline{i}}$.

\begin{proposition}\label{lat_is_stable}
The lattice $Q(\Gamma)_{\pi,\sigma}$ is $W(\Gamma)_{\pi}$-stable.
\end{proposition}
\begin{proof}
It is sufficient to check $\gamma=s_{\overline{i}}(\alpha)\in Q(\Gamma)_{\pi,\sigma}$ for 
any $\alpha \in Q(\Gamma)_{\pi,\sigma}$. If $\pi(i)=i$, then $\gamma_j=\alpha_j, j\neq i$
and $\gamma_i = \alpha_i - 2\langle\alpha ,\varepsilon_i\rangle$.
Since $i$ is adjacent to an even number of vertices,
$\langle\alpha ,\varepsilon_i\rangle$ is integer so $\gamma_i$ is even whenever $\alpha_i$ is.
If $\pi(i)\neq i$, then $s_{\overline{i}}=r_i r_{\pi(i)}$ 
and  $\gamma_j=\alpha_j$ for any $j\neq i,j\neq \pi(i)$.
Since $\alpha$ is symmetric, we have 
\begin{equation}\label{symscal}
\langle \alpha ,\varepsilon_{\pi(i)}\rangle = \langle \alpha ,\varepsilon_i\rangle.
\end{equation}
Note also that $\langle \varepsilon_i, \varepsilon_{\pi(i)}\rangle = 0$, because
$\pi$ is admissible. Hence, 
$\langle r_{\pi(i)}\alpha ,\varepsilon_i\rangle = \langle \alpha ,\varepsilon_i\rangle$,
so $\alpha_i - \gamma_i= \alpha_{\pi(i)} - \gamma_{\pi(i)} = 
2 \langle \alpha ,\varepsilon_i\rangle$.
Thus $\gamma$ is symmetric whenever $\alpha$ is.
\end{proof}

Starting from the graph $\Gamma_{\pi,\sigma}$ one can construct a symmetrizable
GCM matrix $A(\Gamma_{\pi,\sigma})$ (see \cite[$\S$4.7]{kac3}). Namely,
$A(\Gamma_{\pi,\sigma})_{ii}=2$, 
$A(\Gamma_{\pi,\sigma})_{ij}=0$ if $i,j$ are not adjacent in $\Gamma_{\pi,\sigma}$,  
$-1$ for $i - j$ and for $i \Rightarrow j$, 
$-2$ for $i \Leftarrow j$. One can easily check:

\begin{proposition}
$A(\Gamma_{\pi,\sigma})_{ij}=\frac{2\langle \beta_{\overline{i}},\beta_{\overline{j}} \rangle}
{\langle \beta_{\overline{i}},\beta_{\overline{i}} \rangle}.$
\end{proposition}

The root system $\Delta(\Gamma_{\pi,\sigma})$ is defined in terms of the matrix
$A(\Gamma_{\pi,\sigma})$. Namely, we have a lattice $Q(\Gamma_{\pi,\sigma})$
generated by the elements of the simple root system $\Pi(\Gamma_{\pi,\sigma})$ and
$A(\Gamma_{\pi,\sigma})$ gives rise to an integer linear form 
$\langle \beta',\alpha\rangle$ and a reflection 
$s_i'(\alpha) = \alpha - \langle \beta',\alpha\rangle \alpha$
on it for each $\beta'\in \Pi(\Gamma_{\pi,\sigma})$. The Weyl group $W(\Gamma_{\pi,\sigma})$
is generated by these reflections. Comparing the above results, we get:

\begin{corollary}
The natural bijection $\beta_{\overline{i}}\to \beta_i'$ of $\Pi_{\pi,\sigma}$
onto $\Pi(\Gamma_{\pi,\sigma})$ gives rise to an isomorphism of groups
$W(\Gamma)_{\pi}\cong W(\Gamma_{\pi,\sigma})$ and an equivariant
bijection of $(Q(\Gamma)_{\pi,\sigma},W(\Gamma)_{\pi})$ 
onto $(Q(\Gamma_{\pi,\sigma}),W(\Gamma_{\pi,\sigma}))$.
\end{corollary}

Observe that by formula (\ref{symscal}) the condition 
$\langle\alpha ,\varepsilon_i\rangle \leq 0$ for all $i$ is equivalent
for $\alpha\in Q(\Gamma)_{\pi,\sigma}$ to another system of inequalities:
$\langle\alpha ,\beta_{\overline{j}}\rangle \leq 0$ for all $j$. 
Therefore the set $F_{\Gamma,\pi,\sigma}$ of all $\alpha\in Q(\Gamma)_{\pi,\sigma,+}$ 
having a connected support and meeting the above inequalities 
is equal to $F_{\Gamma}\cap \Gamma_{\pi,\sigma}$.
Thus we have:

\begin{corollary}\label{root_in_weyl} The positive root system $\Delta(\Gamma_{\pi,\sigma})_+$ is 
$\Delta(\Gamma_{\pi,\sigma})^{re}_+\cup\Delta(\Gamma_{\pi,\sigma})^{im}_+$, where
\begin{equation}
\Delta(\Gamma_{\pi,\sigma})^{re}_+ = 
   W(\Gamma)_{\pi}\Pi_{\pi,\sigma}\cap Q(\Gamma)_+,
\quad
\Delta(\Gamma_{\pi,\sigma})^{im} = 
   W(\Gamma)_{\pi}F_{\Gamma,\pi,\sigma}.
\end{equation} 
\end{corollary}

We now describe the root system $\Delta(\Gamma_{\pi,\sigma})_+$ in terms of $\Delta(\Gamma)$.
For any $\alpha\in Q(\Gamma)$ set $\overline{\alpha}=\alpha$, if $\alpha\in Q(\Gamma)_{\pi,\sigma}$,
and $\overline{\alpha}=\alpha+\pi(\alpha)$, otherwise. The subsequent Lemma and Proposition
in the case when $\sigma(i)=1$ for all $i\in\Gamma^{\pi}$ follow from \cite[Proposition~4]{hu}.

\begin{lemma}\label{root_system}
$\Delta(\Gamma_{\pi,\sigma})_+=\overline{\Delta(\Gamma)_+} =
\{\overline{\alpha}\vert\alpha\in \Delta(\Gamma)_+\}$.
\end{lemma}
\begin{proof}
We know that $\Delta(\Gamma_{\pi,\sigma})_+=\Delta(\Gamma_{\pi,\sigma})\cap Q(\Gamma)_+$
and analogously $\overline{\Delta(\Gamma)_+}=\overline{\Delta(\Gamma)}\cap Q(\Gamma_+)$.
So it is sufficient to show $\Delta(\Gamma_{\pi,\sigma})=\overline{\Delta(\Gamma)}$.
By Proposition \ref{lat_is_stable} and since the action of $W(\Gamma)_{\pi}$ commutes 
with $\pi$, both sets are $W(\Gamma)_{\pi}$-stable. On the other hand,
$\Pi_{\pi,\sigma}$ and $F_{\Gamma,\pi,\sigma}$ are contained in $\overline{\Delta(\Gamma)_+}$. 
Hence, Corollary  \ref{root_in_weyl} yields the inclusion $\subseteq$ in the Lemma.

Conversely, assume  $\gamma\in\overline{\Delta(\Gamma)_+}$. Apply induction on 
${\rm ht}\gamma= \gamma_1+\cdots+\gamma_n$ to prove that $\gamma$ belongs to 
$\Delta(\Gamma_{\pi,\sigma})_+$. We have three possible cases: 

{\bf 1.} For some $i$, $s_{\overline{i}}(\gamma)\in Q(\Gamma)_+$, 
$\langle \beta_i,\gamma\rangle >0$.
In this case ${\rm ht}s_{\overline{i}}(\gamma)<{\rm ht}\gamma$ and 
$s_{\overline{i}}(\gamma)\in\overline{\Delta(\Gamma)_+}$. Hence, we are done by induction.

{\bf 2.} $\gamma\in F_{\Gamma,\pi,\sigma}$. Clearly, in this case we are done.

{\bf 3.} For some $i$, $s_{\overline{i}}(\gamma)\notin Q(\Gamma)_+$. 
It is well-known that for a root $\alpha$ of $\Gamma$ and a simple reflection
$r_i\in W(\Gamma)$, $r_i\alpha \notin Q(\Gamma)_+$ is equivalent to 
$\alpha=\varepsilon_i$. Applying this observation, we get 
$\gamma=\overline{\varepsilon_i}\in \Delta(\Gamma_{\pi,\sigma})_+$. This completes the proof.
\end{proof}

\begin{proposition}\label{real_1_orbit}
If $\overline{\alpha}\in \Delta(\Gamma_{\pi,\sigma})^{re}_+$, 
then $\alpha\in \Delta(\Gamma)^{re}$ and $\overline{\alpha'}\neq\overline{\alpha}$
for $\alpha'\neq\alpha,\pi(\alpha)$.
\end{proposition}
\begin{proof}
First note that $\Delta(\Gamma_{\pi,\sigma})^{re}$ is divided into 3
disjoint set (possibly empty):
\begin{equation}\label{3_norm}
\Delta(\Gamma_{\pi,\sigma})^{re}=
\bigcup_{\sigma(i)=1}W(\Gamma)_{\pi}\varepsilon_i \sqcup
\bigcup_{j<\pi(j)}W(\Gamma)_{\pi}(\varepsilon_j+\varepsilon_{\pi(j)}) \sqcup
\bigcup_{\sigma(k)=-1}W(\Gamma)_{\pi}2\varepsilon_k.
\end{equation}

\noindent 
Since the action of $W(\Gamma)$ preserves the scalar product we have
$\langle \gamma,\gamma\rangle$ is 1,2,4 for $\gamma$ from the first, second, and third
sets, respectively, because $\langle \varepsilon_j,\varepsilon_{\pi(j)}\rangle=0$.

If $\alpha$ is an imaginary root, then the elements of the smallest height in the
$W(\Gamma)_{\pi}$-orbit of $\overline{\alpha}$ belong to $F_{\Gamma,\pi,\sigma}$
so $\overline{\alpha}$ cannot be real.

Assume that $\gamma=\overline{\alpha_1}=\overline{\alpha_2}$ for two real roots
$\alpha_1$ and $\alpha_2$. If $\gamma$ belongs to the first subset in formula (\ref{3_norm}),
then $\gamma = \alpha_1 = \alpha_2$. Analogously, assuming $\gamma$ to  be in the third
subset, we get $\gamma = 2\alpha_1 = 2\alpha_2$.

Now assume $\gamma=w\varepsilon_j+ w\varepsilon_{\pi(j)}=\alpha+\pi(\alpha)$.
Note that for any $s_{\overline{i}}$, and any $\delta\in \Delta(\Gamma)_+$,
$s_{\overline{i}}(\delta+\pi(\delta))\in Q(\Gamma)_+$ implies 
$s_{\overline{i}}\delta\in Q(\Gamma)_+$. Then $\beta=w^{-1}\alpha\in Q(\Gamma)_+$  
and $\varepsilon_j+ \varepsilon_{\pi(j)}=\beta+\pi(\beta)$.
Thus $\beta$ is either $\varepsilon_j$ or $\varepsilon_{\pi(j)}$.  This completes the proof.
\end{proof}

Assume that a signed quiver
$Q^{\sigma}$ fulfills the condition of \ref{from_sig}; denote by
$\Gamma$ the underlying graph of $\widetilde{Q}$ endowed with an admissible 
involution $\pi=*$. 

\begin{corollary}\label{dim_sym_rep}
Suppose that for any $\alpha\in\Delta(\Gamma)_+\cap Q(\Gamma)_{\pi,\sigma}$ 
either an indecomposable representation of dimension $\alpha$ is symmetric
or there is an indecomposable representation of dimension $\alpha/2$, which
is not isomorphic to any symmetric one. Then:

{\bf 1.} The set of dimensions of indecomposable symmetric representations of $\widetilde{Q}$
is equal to $\Delta(\Gamma_{\pi,\sigma})_+$.

{\bf 2.} If $\gamma\in \Delta(\Gamma_{\pi,\sigma})^{re}_+$, then there exists a unique 
indecomposable symmetric representation of dimension
$\gamma$ modulo isomorphism.
\end{corollary}

\begin{proof}
{\bf 1.} By Proposition \ref{root_system} we may and will replace 
$\Delta(\Gamma_{\pi,\sigma})_+$ by $\overline{\Delta(\Gamma)_+}$.
If $X$ is an indecomposable symmetric representation, then $\dim X$ is sign-matched,
i.e., belongs to $Q(\Gamma)_{\pi,\sigma}$. By Lemma \ref{1_or_2}, either
$X$ is indecomposable as a representation of $\widetilde{Q}$ and in this case
$\dim X=\overline{\dim X}\in \overline{\Delta(\Gamma)_+}$; or $X=V+V^*$ and
$\dim X = \alpha + \pi(\alpha), \alpha = \dim V$. In the latter case, if $\alpha$ 
is not sign-matched,
then $\dim X=\overline{\alpha}\in \overline{\Delta(\Gamma)_+}$. Assume that
$\alpha$ is sign-matched; if $\alpha$ is real, then  $\alpha/2\notin\Delta(\Gamma)$
and by the hypothesis there exists an indecomposable symmetric representation $V'$
of dimension $\alpha$. However, $V'$ must be isomorphic to $V$, since $\alpha$
is a real root. We therefore get a contradiction with  Proposition \ref{2_not_1}.
Thus $\alpha$ is sign-matched and imaginary implying that $\dim X = 2\alpha$ is also
an imaginary root, hence, $\dim X=\overline{\dim X}\in \overline{\Delta(\Gamma)_+}$.

Conversely, if $\gamma\in \overline{\Delta(\Gamma)_+}$, then either 
$\gamma = \alpha +\pi(\alpha)$ with $\alpha$ being a non-sign-matched root
or $\gamma$ itself is a sign-matched root.
In the former case for any indecomposable representation $V$ of dimension $\alpha$,
$V$ is not symmetric, hence, $V+V^*$ is indecomposable symmetric of dimension $\gamma$.
In the latter case the hypothesis yields $\gamma$ is the dimension of an indecomposable
symmetric representation.

{\bf 2.} By Proposition \ref{real_1_orbit}, $\gamma = \overline{\alpha}$ for a unique
(modulo $\pi$) root $\alpha$, and besides $\alpha$ is a real root.
The fact that all indecomposable representations of dimension $\alpha$ are 
conjugate to each other and Corollary \ref{outeriso} complete the proof.
\end{proof}

\section{Signed quivers of finite type.}

\begin{definition}
We call a signed quiver $Q^{\sigma}$ {\it finite} if
there are finitely many isomorphism classes of indecomposable
symmetric representations of $\widetilde{Q}$.
\end{definition}
Note that any usual quiver $S$ gives rise to a signed quiver $S^{\sigma}$
with $\sigma=0$ for all vertices and arrows and no vertices-twins.
Clearly a signed quiver is finite if and only if each connected component
is finite.

\begin{theorem}\label{finite_class}
For a connected signed quiver $Q^{\sigma}$ 3 conditions are equivalent:

(i) $Q^{\sigma}$ is finite

(ii) $\widetilde{Q}$ is finite

(iii) $Q^{\sigma}$ is either a finite quiver or is one of the following signed quivers: 

$$
B_n: {\circ} - \cdots - {\circ}   - \overset{+}{\circ}
\quad\quad\qquad
B_{n,+}: {\circ} - \cdots - {\circ}  \overset{+}{-} {\circ}
$$
 
$$
C_n: {\circ} - \cdots - {\circ}   - \overset{-}{\circ}
\quad\quad\qquad
C_{n,-}: {\circ} - \cdots - {\circ}  \overset{-}{-} {\circ}
$$

\noindent where the orientation of edges is arbitrary, $n$ is the number of vertices.
\end{theorem}

\begin{remark} 
An analogous result for generalized quivers is proved in \cite[3.1]{dw1}.
\end{remark}

\begin{proof}
The implication $(iii)\Rightarrow(ii)$ is obvious. To prove
$(ii)\Rightarrow(iii)$ first observe that if $Q^{\sigma}$ is not a usual quiver, 
then $\widetilde{Q}$ is connected and admits a nontrivial involutive automorphism.
Since $\widetilde{Q}$ is finite, it is of type $A_m,D_m,E_6$. The involutive automorphisms of 
$D_m$ and $E_6$ have a pair of fixed vertices connected by a unique arrow;
this is impossible for $\widetilde{Q}$, because by Definition \ref{sym_quiv} 
we have for this arrow $\varphi$: $t(\varphi^*) = h\varphi, h(\varphi^*) = t\varphi$.
So $\widetilde{Q}$ is $A_m$ with the central symmetry
as $*$ and  $Q$ is $B_{(m+1)/2},C_{(m+1)/2}$ if $m$ is odd and 
$B_{m/2+1,+},C_{m/2+1,-}$ otherwise.

The implication $(ii)\Rightarrow(i)$ follows directly from Lemma \ref{1_or_2}.
Conversely, assume $\widetilde{Q}$ is not finite. Then for an imaginary root
$\alpha$ there exists an infinite family $V(\lambda)$ of pairwise non-isomorphic
indecomposable representations of dimension $\alpha$. Then either of two families, $V(\lambda)$ and
$V(\lambda)+V(\lambda)^*$ contains infinitely many non-isomorphic indecomposable
symmetric representations.
\end{proof}

We want to describe explicitly the indecomposable symmetric representations
of the finite signed quivers. Note that Proposition \ref{embedding}
reduces this question for $B_{n,+}$ to that for $B_n$.
Indeed, we have $B_n = {\rm Ins}_{\varphi}(B_{n,+})$, where $\varphi$ denote the 
signed arrow of $B_{n,+}$. Analogously, $C_n = {\rm Ins}_{\varphi}(C_{n,-})$
and it is sufficient to consider  $B_n$ and $C_n$.

By Lemma \ref{1_or_2}, to describe the indecomposable symmetric representations
of $Q^{\sigma}$ we only need to answer a question as follows: which
indecomposable representations $V$ of $\widetilde{Q}$ admit a signed form
making $V$ symmetric? Clearly, a necessary condition is that $\dim V$ is sign-matched.
In both cases  $\widetilde{Q}=A_{2n+1}$ and the indecomposable representations
have dimensions of form $(0,\cdots,0,1,\cdots,1,0,\cdots,0)$. So for $C_n$ none
root is sign-matched, because the middle dimension is odd. For $B_n$ the symmetric
dimensions are sign-matched and clearly, the corresponding indecomposable representations
are symmetric with respect to an appropriate signed form. Note that for both $B_n$ and $C_n$
the involution $*$ is admissible and the above description shows that the hypothesis
of Corollary \ref{dim_sym_rep} holds. Note also that the graph $\Gamma_{\pi,\sigma}$ is
the Dynkin diagram $B_n$ and $C_n$ for the quiver $B_n$ and $C_n$, respectively.  
Hence, we obtain:

\begin{theorem}\label{finite_dim}
The set of dimensions of indecomposable symmetric representations
of $B_n$ and $C_n$ constitute a root system of type
$B_n$ and $C_n$, respectively. 
\end{theorem}

\begin{remark} 
An analogous result for generalized quivers is proved in \cite[3.6]{dw1}.
\end{remark}

\section{Signed quivers of tame type.}

\begin{definition}
We call a signed quiver $Q^{\sigma}$ {\it tame} if
infinite families of classes of indecomposable
symmetric representations of $\widetilde{Q}$ exist and depend on 1 parameter.
\end{definition}

\begin{theorem}\label{tame_class}
For a connected signed quiver $Q^{\sigma}$ 3 conditions are equivalent:

(i) $Q^{\sigma}$ is tame

(ii) $\widetilde{Q}$ is tame

(iii) $Q^{\sigma}$ is either a tame quiver or is one of the following  signed quivers: 

\begin{equation}\label{loops}
O_+:\overset{\circ+}{\circlearrowleft +} \qquad  O_-: \overset{\circ+}{\circlearrowleft -} \qquad
Sp_-: \overset{\circ-}{\circlearrowleft -} \qquad Sp_+:\overset{\circ-}{\circlearrowleft +} 
\end{equation}

$$
D_n^{(2)}: \overset{+}{\circ} -  {\circ} - \cdots - {\circ}   - \overset{+}{\circ}
\quad
D_{n,+}^{(2)}: {\circ} \overset{+}{-} {\circ} - \cdots - {\circ}   - \overset{+}{\circ}
\quad
D_{n,++}^{(2)}: {\circ} \overset{+}{-}  {\circ} - \cdots - {\circ} \overset{+}{-}  {\circ}
$$

$$
C_n^{(1)}: \overset{-}{\circ} - {\circ} - \cdots - {\circ}   - \overset{-}{\circ}
\quad
C_{n,-}^{(1)}: {\circ}\overset{-}{-} {\circ} - \cdots - {\circ}   - \overset{-}{\circ}
\quad
C_{n,--}^{(1)}: {\circ} \overset{-}{-} {\circ} - \cdots - {\circ} \overset{-}{-} {\circ}
$$

$$
A_{2n}^{(2)}: \overset{-}{\circ} - {\circ} - \cdots - {\circ}   - \overset{+}{\circ}
$$

$$
A_{2n,+}^{(2)}: \overset{-}{\circ} - {\circ} - \cdots - {\circ} \overset{+}{-} {\circ}
\quad
A_{2n,-}^{(2)}: {\circ} \overset{-}{-}  {\circ} - \cdots - {\circ}   - \overset{+}{\circ}
\quad
A_{2n,+-}^{(2)}: {\circ} \overset{-}{-}  {\circ} - \cdots - {\circ}  \overset{+}{-} {\circ}
$$

$$
Z_n: {\circ}_1 - {\circ} - \cdots - {\circ}   - {\circ}_{1^*}
$$

\begin{alignat*}{2}
 B_n^{(1)}:  \overset{+}{\circ} - {\circ} - \cdots - & {\circ} - {\circ}  &
 \qquad  B_{n,+}^{(1)}:   {\circ}\overset{+}{-} {\circ} - \cdots - & {\circ} - {\circ}
\\
    & \,\vert  &\qquad  & \,\vert\\
    & {\circ} & \qquad    & {\circ} 
\end{alignat*}

\begin{alignat*}{2}
 A_{2n-1}^{(2)}:  \overset{-}{\circ} - {\circ} - \cdots - & {\circ} - {\circ}  &
 \qquad A_{2n-1,-}^{(2)}:  {\circ}\overset{-}{-} {\circ} - \cdots - & {\circ} - {\circ}
\\
    & \,\vert  &\qquad  & \,\vert\\
    & {\circ} & \qquad    & {\circ} 
\end{alignat*}

\noindent where the orientation of edges is arbitrary, the number of vertices is $n+1$.
\end{theorem}
\begin{remark}
A similar result for generalized quivers is proved in \cite[4.2]{dw1}. However,
the signed
quivers $O_+,Sp_-,D_{n,+}^{(2)},C_{n,-}^{(1)},A_{2n}^{(2)},A_{2n,+-}^{(2)}$ have no
counterpart in that theory (see Remark \ref{DW-gq}).
\end{remark}
\begin{proof}
We argue exactly as in the proof of \ref{finite_class}. In the proof of implication
$(ii)\Rightarrow(iii)$ we have 4 possibilities for the underlying graph $\Gamma$ 
of $\widetilde{Q}$ and the involutive authomorphism $\pi$ of $\Gamma$:

 $(a):\Gamma$ is a graph with 1 vertex and one loop, $\pi$ is trivial;

 $(b):\Gamma$ is a circular graph $A_m^{(1)}$ ($m+1$ vertices), $\pi$ is a mirror symmetry;

 $(c):\Gamma$ is a circular graph $A_m^{(1)}$, $m$ is odd, 
 $\pi$ is the central symmetry;
  
 $(d):\Gamma$ is $D_n^{(1)}$, $\pi$ is the central symmetry.

\noindent The case $(a)$ yields the signed quivers from (\ref{loops}). In the case
$(b)$ there are several possibilities for the mirror: if the mirror meets
two vertices, we get $D_n^{(2)},C_n^{(1)},A_{2n}^{(2)}$; if the mirror meets
a vertex and the opposite edge, we get 
$D_{n,+}^{(2)},C_{n,-}^{(1)},A_{2n,+}^{(2)},A_{2n,-}^{(2)}$;
finally, when the mirror meets two edges, we get $D_{n,++}^{(2)},C_{n,--}^{(1)},A_{2n,+-}^{(2)}$.
The case $(c)$ yields $Z_{(m+1)/2}$. The case $(d)$ with center being a vertex yields
$B_n^{(1)},A_{2n-1}^{(2)}$ and with center being on an edge yields
$B_{n,+}^{(1)},A_{2n-1,-}^{(2)}$.
\end{proof}

We are now going to describe the indecomposable symmetric representations
for the tame signed quivers. Note that the notation introduced
in \ref{tame_class} divide them into families with
a base term $X$ and  additional terms $X_+,X_-$ etc.
For example, $A_{2n}^{(2)}$  has 3 relatives: 
$A_{2n,+}^{(2)},A_{2n,-}^{(2)},A_{2n,+-}^{(2)}$.
These quivers can be obtained from  each other applying repeatedly the functors Del
and Ins: 
${\rm Ins}_{\varphi}(A_{2n,+-}^{(2)}) = A_{2n,-}^{(2)}$,
${\rm Ins}_{\psi}(A_{2n,+-}^{(2)}) = A_{2n,+}^{(2)}$,
${\rm Ins}_{\varphi}(A_{2n,+}^{(2)}) = A_{2n}^{(2)}=
{\rm Ins}_{\psi}(A_{2n,-}^{(2)})$, where $\varphi$ and $\psi$ are the + and the - arrow
of $A_{2n,+-}^{(2)}$, respectively.
So if we know the indecomposable symmetric representations for $A_{2n}^{(2)}$,
Proposition \ref{embedding} yields the same for other 3 quivers.
We therefore need to consider the quivers from (\ref{loops}) and 
$D_n^{(2)},C_n^{(1)},A_{2n}^{(2)},Z_n,B_n^{(1)},A_{2n-1}^{(2)}$.

By Lemma \ref{1_or_2}, to describe the indecomposable symmetric representations
of $Q^{\sigma}$ is equivalent to find the indecomposable representations $V$ 
of $\widetilde{Q}$ such that $\dim V$ is sign-matched and
$V$ is symmetric with respect to an appropriate signed form.
First we consider the quivers from (\ref{loops}). Here $\widetilde{Q}$ is the
quiver with one arrow-loop, an indecomposable representation $A$
is just a Jordan matrix $J_n(\lambda)$.

\begin{example} $Q^{\sigma} = O_+$. A signed form for $A$
is a symmetric matrix $J$ such that $JA$ is symmetric (cf. \ref{++}). 

\begin{proposition}\label{symsym}
Any matrix $A$ has a presentation $JA=B$ with $J^{\top}=J$, $\det J\neq 0$,
$B^{\top}=B$.
\end{proposition}

\begin{proof}
We can rewrite $A=CB$ with $C=J^{-1}, C^{\top} = C$. 
Furthermore, a standard fact from Linear Algebra is that
$C$ can be presented as $C=D D^{\top}$, hence,
$A=D D^{\top} B  D D^{-1}$. Our presentation 
is therefore equivalent to the classical fact that $A$ is conjugate to a symmetric matrix
$D^{\top} B  D$ (see. e.g. \cite{ga}). 
\end{proof}
\end{example}

\begin{example} $Q^{\sigma} = O_-$. A signed form for $A$
is a symmetric matrix $J$ such that $JA$ is anti-symmetric (cf. \ref{+-}). 
Arguing as in the proof of \ref{symsym}, we get: 

\begin{proposition}\label{symskew}
Jordan matrix $A=J_n(\lambda)$ has a presentation $JA=B$ with $J^{\top}=J,\det J\neq 0$
$B^{\top}=-B$ if and only if $\lambda=0$ and $n$ is odd.
\end{proposition}
\end{example}

\begin{example} $Q^{\sigma} = Sp_+$. We consider only representations
of sign-matched dimensions, so $n$ is even.  A signed form for $A$
is an anti-symmetric matrix $J$ such that $JA$ is symmetric. 

\begin{proposition}\label{skewsym}
Jordan matrix $A=J_n(\lambda)$ has a presentation 
$JA=B$ with $J^{\top}=-J,\det J\neq 0$
$B^{\top}=B$ if and only if $\lambda=0$ and $n$ is even.
\end{proposition}
\begin{proof}
The property is equivalent to $A=CB$ with $C^{\top} = -C, \det C\neq 0$.
Using the properties of trace, we get: $tr(A)=tr(CB)^{\top}=-tr(BC)=-tr(A)$.
Hence, $tr(A)=n\lambda=0$. In addition, $\det C\neq 0$ implies
$n$ is even. Conversely, given an even $n$ we take $C'$ and $B'$ such that
$C_{ij}'=0$ for $i+j\neq n+1$ and  $B_{ij}'=0$ for $i+j\neq n+2$. Then $C'B'$
is conjugate to $J_n(0)$, hence, $J_n(0)$ can also be presented in this form.
\end{proof}
\end{example}

\begin{example} $Q^{\sigma} = Sp_-$. Again $n$ is even.  
A signed form for $A$ is an anti-symmetric matrix $J$ such that $JA$ is anti-symmetric. 

\begin{proposition}\label{skewskew}
None Jordan matrix $A=J_n(\lambda)$ with even $n$ has a presentation 
$JA=B$ with $J^{\top}=-J,\det J\neq 0$
$B^{\top}=-B$.
\end{proposition}

\begin{proof}
We have: $B=J(\lambda{\rm Id}+J_n(0))$ and $B-\lambda J$ is anti-symmetric,
so we may assume $\lambda=0, B=J J_n(0)$. However, $i$-th column of 
$J_n(0)$ is $i-1$-th column of $J$ for $i>1$ and the first one is zero.
Using the condition that both $J$ and $B$ are anti-symmetric, we get
$B=J=0$. 
\end{proof}
\end{example}

Now we consider $D_n^{(2)},C_n^{(1)},A_{2n}^{(2)},B_n^{(1)},A_{2n-1}^{(2)}$.
Observe that all these quivers with $n\geq 2$ fulfill the condition from
\ref{from_sig} and the involution $\pi$ of the underlying graph
$\Gamma$ of $\widetilde{Q}$ is admissible. As diagram
$\Gamma_{\pi,\sigma}$ we obtain the affine graph denoted by exactly
the same term as the quiver: $D_n^{(2)}$ for $D_n^{(2)}$ etc.
(it is of course more correct to say that we denoted
our signed quivers following the notation of Kac, see \cite{kac2},\cite{kac3}).
For the signed quiver $Z_n,n\geq 2$ the involution is also admissible and
we have  $\Gamma_{\pi,\sigma}=A_{n-1}^{(1)}$.
Below we consider these quivers case by case and prove 

\begin{theorem}\label{tame_dim}
Let $Q^{\sigma}$ be one
of the quivers $D_n^{(2)},C_n^{(1)},A_{2n}^{(2)},Z_n,B_n^{(1)},A_{2n-1}^{(2)}$,
$n\geq 2$.
Then the set of dimensions for the indecomposable symmetric
representations of $\widetilde{Q}$ is equal to the root system 
$\Delta(\Gamma_{\pi,\sigma})_+$
of the corresponding affine graph $\Gamma_{\pi,\sigma}$. Furthermore,
all indecomposable symmetric representations with dimension
from $\Delta(\Gamma_{\pi,\sigma})_+^{re}$ are isomorphic to each other.
\end{theorem}
\begin{proof}
For almost all cases the Theorem follows from Corollary \ref{dim_sym_rep}
and 

\begin{proposition}\label{5diagram}
The hypothesis of \ref{dim_sym_rep} holds for the
quivers $D_n^{(2)},C_n^{(1)}, A_{2n}^{(2)}$, $B_n^{(1)},A_{2n-1}^{(2)}, n\geq 2$, i.e.,
for any sign-matched root $\alpha$
either an indecomposable representation of dimension $\alpha$ is symmetric
or there is an indecomposable representation of dimension $\alpha/2$, which
is not isomorphic to any symmetric one. 
\end{proposition}

\begin{proof} {\it (of Proposition )} 
For $Q^{\sigma} = D_n^{(2)},C_n^{(1)},A_{2n}^{(2)}$ the underlying graph of 
$\widetilde{Q}$ is $A_{2n-1}^{(1)}$ such that the vertices in the clockwise order are 
$1,2,\cdots,n+1,n^*,\cdots,2^*$.
For any $i\in \widetilde{Q}_0$, write $i+1$ and $i-1$ for the next and the
previous vertex in the clockwise order, respectively;
denote by $\varphi_i\in \widetilde{Q}_1$
the unique arrow such that $\{t\varphi_i, h\varphi_i\}=\{i,i+1\}$.
We keep notation $J_d(\lambda)$ for the corresponding Jordan matrix;
by $I_{r,s}$ we denote the $s\times r$ matrix having
the first $p=min(r,s)$ standard basis vectors of ${\bf k}^s$ 
as the first $p$ columns, the last $r-p$ columns being zero; 
by $I_{r,s}'$ we denote the $s\times r$ matrix having
the {\it last} $p$ standard basis vectors of ${\bf k}^s$ as
the last $p$ columns, the first $r-p$ columns being zero; 
for the identity matrix $I_{d,d}$ we just write $I_d$.

The indecomposable representations of $\widetilde{Q}$ are as follows. 
The root $\delta = (1,\cdots,1)$ generates
the fundamental domain $F_{\widetilde{Q}}$. For any natural $d$ we have a
1-dimensional family $V_{\lambda}^d,\lambda\in {\bf k}$ of indecomposable 
representations of dimension $d\delta$:
\begin{equation}\label{A_m-imag}
V_{\lambda}^d(\varphi_1) = J_d(\lambda), V_{\lambda}^d(\varphi_i) = I_d, i\neq 1.
\end{equation}
The real roots are the vectors $\alpha^{k,l,d}$ such that: 
$(\alpha^{k,l,d}- d\delta)_i=1$, for
$i=k,k+1,\cdots,l$ and 0 for other vertices, $k,l\in \widetilde{Q}_0$,
$d\in {\bf Z}_+$. For each $\alpha^{k,l,d}$ we
have an indecomposable representation $V^{k,l,d}$ as follows:
\begin{equation}\label{A_m-real}
V^{k,l,d}(\varphi)=I_{r,s},\varphi\in\widetilde{Q}_1\setminus\{\varphi_l\},
V^{k,l,d}(\varphi_l)=I_{r,s}',\; r = \alpha^{k,l,d}_{t\varphi},
s=\alpha^{k,l,d}_{h\varphi}.
\end{equation}
Note that $V^{k,l,d}(\varphi_i)$ is the identity matrix for all $i$ but
$l$ and $k-1$.

{\bf 1.} $Q^{\sigma} = D_n^{(2)}$. The sign-matched dimension vectors are just
symmetric ones; so the sign matched roots are $d\delta$ and $\alpha^{i,i^*,d}$
for arbitrary $i$ and $d$. We claim that the corresponding representations
$V_{\lambda}^d$ and $V^{i,i^*,d}$ are symmetric with respect to an appropriate
signed form.

Let $V=V_{\lambda}^d$. Then in the dual basises of $V(i)^*$ the representation
$V^*$ is given by the matrices as follows: $V^*(\varphi)=I_d$ for $\varphi\neq\varphi_{2^*}$,
and $V^*(\varphi_{2^*})=J_d^{\top}(\lambda)$. We are looking for a signed form
$J_k:V(k)\to V(k^*)^*$ such that $J_1^*=J_1,J_{n+1}^*=J_{n+1}$, $J_i^*=J_{i^*}$ for
$i=2,\cdots,n$, and $J$ is an isomorphism of $V$ onto $V^*$.
Taking into account the explicit form of $V$ and $V^*$ we rewrite the condition:
$J_2=J_3=\cdots=J_{n+1}=J_{n^*}=\cdots=J_{2^*}$, $J_1^{\top}=J_1,J_2^{\top}=J_2$, 
and the diagram commutes:
$$
\begin{CD}
V(2^*)         @>I_d>>   V(1)    @>J_d(\lambda)>>  V(2) \\
@VJ_2^{\top}VV  @VVJ_1V    @VVJ_2V\\
V(2)^*     @>>J_d^{\top}(\lambda)> V(1)^*  @>>I_d> V(2^*)^*
\end{CD}
$$
\noindent The left square condition is equivalent to the right one:
$J_2 J_d(\lambda)= J_1$. By Proposition \ref{symsym},
such $J_1$ and $J_2$ exist.

Let $V=V^{i,i^*,d},\alpha=\alpha^{i,i^*,d}$. Note that $(i-1)^*=i^*+1$ 
and $\varphi_{i-1}^*=\varphi_{i^*}$.
So $V(\varphi)$ is the identity matrix for all arrows but $\psi=\varphi_{i-1}$
and $\psi^*$, $V(\psi)=I_{r,s},V(\psi^*)=I_{s,r}'$, where 
$r=\alpha_{t\psi}=\alpha_{h\psi^*}, s=\alpha_{h\psi}=\alpha_{t\psi^*}$.
Hence, $V^*(\varphi)=I_d$ or $I_{d+1}$ for $\varphi\neq\psi,\psi^*$, 
$V^*(\psi)=I_{r,s}',V^*(\psi^*)=I_{s,r}$. So the condition on the signed
form $J_k$ is equivalent to finding a symmetric $d\times d$ matrix $J_1$
and symmetric $d+1\times d+1$ matrix $J_2$ such that $I_{d,d+1}'J_1=J_2I_{d,d+1}$.
We may take both $J_1$ and $J_2$ having units on the secondary diagonal and 0 outside it.
Thus the Proposition is proved for this case.

{\bf 2.} $Q^{\sigma} = C_n^{(1)}$. A sign-matched dimension vector $\alpha$ is a
symmetric one such that $\alpha_1$ and $\alpha_{n+1}$ are both even.
So the only sign-matched roots are $2m\delta$ for arbitrary $m$. 
We claim that none representation $V_{\lambda}^d$ is symmetric; clearly,
this implies what we need. Assume that  $d$ is even, 
otherwise the claim is obvious. For $V=V_{\lambda}^d$, $V^*$ is as in case {\bf 1} but
for two arrows with heads $1,n+1$ the maps are multiplied by -1
(see Definition \ref{v^*}). Simlarly to case {\bf 1}, a signed form
can be reduced to anti-symmetric non-degenerate $d\times d$ matrices $J_1,J_2$ such that
$J_2 J_d(\lambda) = J_1$. By Proposition \ref{skewskew}, this is impossible.

{\bf 3.} $Q^{\sigma} = A_{2n}^{(2)}$. A sign-matched dimension vector $\alpha$ is a
symmetric one such that $\alpha_1$ is even. So we have more possibilities for
$\alpha$: $2m\delta$ and $\alpha^{i,i^*,d}$, where either $d$ is even and $i=2,\cdots,n+1$,
or $d$ is odd and $i=n^*,\cdots, 2^*,1$. We claim that there are no symmetric
representations isomorphic to $V_{\lambda}^{2m}$, whereas for the real sign-matched
roots a symmetric representation exists.

Let $V=V_{\lambda}^{2m}$. Then, as in cases {\bf 1,2} a signed form
can be reduced to non-degenerate $d\times d$ matrices $J_1,J_2$, such that $J_1^{\top}=-J_1$,
$J_2^{\top}=J_2$, and  $J_2 J_{2m}(\lambda) = J_1$. By Proposition \ref{symskew}, 
this is impossible.

Let $V=V^{i,i^*,d},\alpha=\alpha^{i,i^*,d}$, where $d$ is even and $i=2,\cdots,n+1$. 
As in case {\bf 1} a signed form can be reduced to an anti-symmetric 
$d\times d$ matrix $J_1$ and a symmetric $d+1\times d+1$ matrix $J_2$ such that 
$I_{d,d+1}'J_1=J_2I_{d,d+1}$. Take $J_1$ having a sequence $1,-1,1,\cdots,-1$ 
on the secondary diagonal and 0 outside it; as $J_2$ take the submatrix of $J_1$ obtained
removing first row and last column. The case when $d$ is odd and $i=n^*,\cdots, 1^*,1$ is similar.

For the quivers $Q^{\sigma}=B_n^{(1)},A_{2n-1}^{(2)}$ the underlying graph of 
$\widetilde{Q}$ is $D_{2n}^{(1)}$ with vertices $1,2,3,\cdots,n,n+1,n^*,\cdots,3^*,2^*,1^*$.
Assume that $1,2$ are tails of the arrows $\chi,\psi\in\widetilde{Q}_1$ respectively
such that $h\chi=h\psi=3$ (this is possible not for each orientation
but the consideration of other cases is similar).

The symmetric roots are as follows. The root $\delta=(1,1,2,\cdots,2,1,1)$ generates
the fundamental domain $F_{\widetilde{Q}}$. The real symmetric roots are of
form $d\delta+\alpha$, where either $d\geq 0$, each $\alpha_i$ is 0 or 1 
and the support of $\alpha$ is a connected symmetric subgraph in  $D_{2n}^{(1)}$, 
or $d>0$ and $\alpha=-\varepsilon_i - \cdots - \varepsilon_{i^*}$ with $i=4,\cdots,n+1$.
In dimension $d\delta$ we construct a 1-parameter family $V_{\lambda}^d$ of 
indecomposable representations:
\begin{equation}\label{D_n-imag}
V_{\lambda}^d(\chi)=(I_d\vert 0)^{\top},V_{\lambda}^d(\psi)=(0\vert I_d)^{\top},
V_{\lambda}^d(\chi^*)=(I_d\vert I_d),V_{\lambda}^d(\psi^*)=(I_d\vert J_d(\lambda)),
\end{equation}
\noindent $V_{\lambda}^d(\varphi)=I_{2d}$ for all other arrows $\varphi$.

{\bf 4.} Let $Q^{\sigma}=B_n^{(1)}$. Here all the above symmetric roots are sign-matched.
We claim that all the corresponding indecomposable representations are symmetric
with respect to an appropriate signed form.

Let $V=V_{\lambda}^d$. Then in the dual basises of $V(i)^*$ the representation
$V^*$ is: 
\begin{equation}\label{D_n-imag^*}
V^*(\chi)=(I_d\vert I_d)^{\top},V^*(\psi)=(I_d\vert J_d(\lambda))^{\top},
V^*(\chi^*)=(I_d\vert 0),V^*(\psi^*)=(0\vert I_d),
\end{equation}
\noindent $V_{\lambda}^d(\varphi)=I_{2d}$ for all other arrows $\varphi$.
For a signed form $J\in{\rm Hom}(V,V^*)$ we therefore have:
$J_3=J_4=\cdots = J_{3^*}=A$, where $A$ is a symmetric $2d\times 2d$
matrix. Denote by $A_{11}=A_{11}^{\top}, A_{12}, A_{21}=A_{12}^{\top}, A_{22}=A_{22}^{\top}$
the $d\times d$ blocks of $A$. Then a direct calculation
shows that $J$ is a signed form for $V$ if and only if:
$J_1=J_2=A_{11}=A_{12}=A_{21}$ and $A_{22}=J_d(\lambda) A_{11}$. By Proposition \ref{symsym},
we may find such symmetric matrices $A_{11}$ and $A_{22}$. 

As for the real roots, we will not present calculations for all of them,
because the notation needs much place. We just consider the case
$n=2,\alpha=(0,0,1,0,0)$ and we assure a reader that the other
cases are not too different. As an indecomposable representation
of dimension $d,d,2d+1,d,d$ we can take $V$ such that:
\begin{equation}\label{D_4^1-real}
V(\chi)=(I_d\vert 0)^{\top},V(\psi)=(I_d\vert I_{d+1,d})^{\top},
V(\chi^*)=(I_d\vert 0),V(\psi^*)=(I_d\vert I_{d+1,d}').
\end{equation}
\noindent Let a collection $(J_1,J_2,A=J_3,J_1^{\top},J_2^{\top})$ be a candidate
for a signed form. Decompose $A$ into blocks with respect to the
decomposition of ${\bf k}^{2d+1}={\bf k}^d\oplus {\bf k}^{d+1}$ used
in (\ref{D_4^1-real}): $A_{11}=A_{11}^{\top}, A_{12}, A_{21}=A_{12}^{\top}, A_{22}=A_{22}^{\top}$.
Then a direct calculation shows that $J$ is a signed form for $V$ if and only if:
$J_1=J_2=A_{11}$, $A_{12}=A_{21}=0$, and $I_{d,d+1}'A_{11}=A_{22}I_{d,d+1}$.
We indicated the appropriate $A_{11}$ and $A_{22}$ in case {\bf 1}.

{\bf 5.} Let $Q^{\sigma}=A_{2n-1}^{(2)}$. Here a symmetric vector $\alpha$ is sign-matched
if and only if $\alpha_{n+1}$ is even. Hence, the only sign-matched roots are the
multiples of $\delta$.
We claim that none representation $V=V_{\lambda}^d$ is symmetric
with respect to a signed form. Formula (\ref{D_n-imag^*})
holds for $V^*$ and for all other arrows $V^*(\varphi)=I_{2d}$ except
for the arrow $\rho$ having $n+1$ as the head, where $V^*(\rho)=-I_{2d}$.
Arguing as in {\bf 4} we reduce a signed form
to a presentation $A_{22}=J_d(\lambda) A_{11}$, where
$A_{11}$ and $A_{22}$ are anti-symmetric non-degenerate
matrices. By Proposition \ref{skewskew}, such a presentation does not exist. 

Thus we considered all 5 signed quivers and the Proposition is proved.
\end{proof}

To finish the proof of the Theorem we need to consider the quiver $Q^{\sigma} = Z_n$. 
In this case, $\widetilde{Q}= A_{2n-1}^{(1)}$ such that the vertices in the 
clockwise order are $1,2,\cdots,n,1^*$, $2^*,\cdots,n^*$. 
The involution $\pi$ takes each vertex to its opposite and 
each symmetric dimension vector is sign-matched.
The roots of type $\alpha^{k,l,d}$ are therefore not symmetric, 
so the only sign-matched roots are the multiples of $\delta$.

Let $V=V_{\lambda}^d$. We claim: $V$ is symmetric with respect to
a signed form if and only if either $\lambda=1$ and $d$ is odd or
$\lambda=-1$ and $d$ is even. Actually, arguing as in the proof
of \ref{5diagram}, we see that a signed form can be reduced
to a non-degenerate matrix $J$ such that $J J_d(\lambda) = J^{\top}$.
Set $A=J_d(\lambda)$ and rewrite the condition as $A = J^{-1}J^{\top}$;
then we get $A^{-\top}= J^{\top}J^{-1}$, hence, $tr(A)=tr(A^{-\top})$ and $\lambda = \lambda^{-1}$.
So we have: 
\begin{equation}\label{strange}
\lambda J+J J_d(0)=J^{\top},\lambda=\pm 1,\det J\neq 0.
\end{equation}
If $\lambda=1$, this implies $d-1={\rm rank} J J_d(0)$ is even, since $J^{\top} - J$ is anti-symmetric.
If $\lambda=-1$, one can show that any solution of (\ref{strange}) has the
properties: $J_{ij}=0$, if $i+j\leq d$, $J_{ij}=-J_{ji}$, if $i+j=d+1$; so $\det J\neq 0$ 
implies $d$ is even. A reader can check that under the above restrictions (\ref{strange}) has
a solution.

Now we see that the assertion of \ref{dim_sym_rep} is not true only for the imaginary
roots $\gamma=(2m+1)\delta$. For this case we can however present $\gamma$ as
the symmetrization $\alpha+\overline{\alpha}$ of a real root $\alpha$ such that
$\alpha_i=m+1$, $\alpha_{j^*}=m$. Thus the proof of \ref{dim_sym_rep} yields
the assertion of the Theorem for $Z_n$.
\end{proof}


\end{document}